\documentclass{amsart}
\usepackage{amsfonts,amssymb,epsfig}
\usepackage[latin1]{inputenc}
\usepackage[all]{xy}
\usepackage{amsmath}
\usepackage{delarray}
\usepackage{amssymb}
\newtheorem{thm}{Theorem}[section]
\newtheorem{cor}[thm]{Corollary}
\newtheorem{lem}[thm]{Lemma}
\newtheorem{prop}[thm]{Proposition}
\theoremstyle{definition}
\newtheorem{defn}[thm]{Definition}
\theoremstyle{remark}
\newtheorem{rem}[thm]{Remark}
\numberwithin{equation}{section}
\newtheorem{example}[thm]{Example}

\renewcommand{\epsilon}{\varepsilon}

\newcommand{\fg}{\mathfrak{g}}
\newcommand{\fh}{\mathfrak{h}}

\newcommand{\cL}{\mathcal L}

\newcommand{\N}{{\mathbb N}}
\newcommand{\Z}{{\mathbb Z}}
\newcommand{\K}{{\mathbb K}}
\newcommand{\R}{{\mathbb R}}
\newcommand{\C}{{\mathbb C}}

\newcommand{\QED}{\hfill $\square$\vspace{2mm}}

\begin{document}

\title[Normal forms of vector fields on Poisson manifolds]
{Normal forms of vector fields on Poisson manifolds}

\author{Philippe Monnier}
\address{Laboratoire Emile Picard, UMR 5580 CNRS, UFR MIG, Université
Toulouse III} \email{pmonnier@picard.ups-tlse.fr}
\author{Nguyen Tien Zung}
\address{Laboratoire Emile Picard, UMR 5580 CNRS, UFR MIG, Université
Toulouse III} \email{tienzung@picard.ups-tlse.fr}


\subjclass{??}

\date{July/2005}

\begin{abstract}
We study formal and analytic normal forms of radial and Hamiltonian vector
fields on Poisson manifolds near a singular point.
\end{abstract}

\maketitle
\section{Introduction}

This paper is devoted to the study of normal forms à la Poincaré-Birkhoff
for analytic or formal vector fields on Poisson manifolds. We will be
interested in two kinds of vector fields, namely Hamiltonian vector
fields, and ``radial'' vector fields, i.e. those vector fields $X$ such
that $[X,\Pi] = \cL_X\Pi = - \Pi$, where $\Pi$ denotes the Poisson
structure, and the bracket is the Schouten bracket. Our motivation for
studying radial vector fields comes from Jacobi structures
\cite{Marle-Jacobi}, while of course the main motivation for studying
Hamiltonian vector fields comes from Hamiltonian dynamics. We will assume
that our vector field $X$ vanishes at a point, $X(0) = 0$, and that the
linear part of $\Pi$ or of its transverse structure at $0$ corresponds to
a semisimple Lie algebra. In this case, it is well known
\cite{Weinstein-Poisson1983,Conn-Analytic1984} that $\Pi$ admits a formal
or analytic linearization in a neighborhood of $0$. We are interested in a
simultaneous linearization or normalization of $\Pi$ and $X$.

In Section \ref{section:Poisson-homo}, we study the problem of
simultaneous linearization of couples $(\Pi,X)$ where $\Pi$ is a Poisson
structure and $X$ is a vector field such that $\cL_X\Pi = - \Pi$. Such
couples are called \emph{homogeneous Poisson structures} in the sense of
Dazord, Lichnerowicz and Marle \cite{Marle-Jacobi}, and they are closely
related to Jacobi manifolds. Such structures are closely related to Jacobi
structures. More precisely, a 1-codimensional submanifold of a homogeneous
Poisson manifold $(M,\Pi,X)$ which is transverse to the vector field $X$
has an induced Jacobi structure, and all Jacobi manifolds can be obtained
in this way. On the other hand, a 1-codimensional submanifold of a Jacobi
manifold $(N,\Lambda,E)$ transverse to the structural vector field $E$ has
an induced homogeneous Poisson structure, and all homogeneous Poisson
manifolds can be obtained in this way (see \cite{Marle-Jacobi}). Our first
result is the following (see Theorem \ref{thm:Poisson-hom-formal}):

{\bf Theorem A}. {\it Let $(\Pi,X)$ be a formal homogeneous Poisson
structure on $\K^n$ (where $\K$ is $\C$ or $\R$) such that the linear part
$\Pi_1$ of $\Pi$ corresponds to a semisimple Lie algebra $\fg$. Suppose
that its linear part $(\Pi_1,X^{(1)})$ is semisimple nonresonant. Then
there exists a formal diffeomorphism which sends $(\Pi,X)$ to
$(\Pi_1,X^{(1)})$.}

The \emph{semisimple nonresonant} condition in the above theorem is a
generic position on $X^{(1)}$: the set of $X^{(1)}$ which does not satisfy
this condition is of codimension 1, and moreover if $X^{(1)} - I$ is
diagonalizable and small enough, where $I = \sum x_i {\partial \over
\partial x_i}$ denotes the standard radial (Euler) vector field, then the
semisimple nonresonant is automatically satisfied.

For analytic linearization, due to possible presence of small divisors, we
need a Diophantine-type condition. Here we choose to work with a modified
Bruno's $\omega$-condition \cite{Bruno,Bruno2} adapted to our case. See
Definition \ref{defn:omega_condition} for the precise definition of our
$\omega$-condition. The set of $(\Pi_1,X^{(1)})$ which satisfy this
$\omega$-condition is of full measure. We have (see Theorem
\ref{thm:Poisson-hom-analytic}):

{\bf Theorem B}. {\it Let $(\Pi,X)$ be an analytic homogeneous Poisson
structure on $\K^n$ (where $\K$ is $\C$ or $\R$) such that the linear part
$\Pi_1$ of $\Pi$ corresponds to a semisimple Lie algebra. Suppose moreover
that its linear part $(\Pi_1,X^{(1)})$ is semisimple nonresonant and
satisfies the $\omega$-condition. Then there exists a local analytic
diffeomorphism which sends $(\Pi,X)$ to $(\Pi_1,X^{(1)})$.}

In Section \ref{section:Hamiltonian}, we study local normal forms of
Hamiltonian systems on Poisson manifolds. According to Weinstein's
splitting theorem \cite{Weinstein-Poisson1983}, our local Poisson manifold
$((\K^n,0),\Pi)$, where $\K = \R$ or $\C$, is a direct product
$((\K^{2l},0),\Pi_{symp}) \times ((\K^m,0),\Pi_{trans})$ of two Poisson
manifolds, where the Poisson structure $\Pi_{symp}$ is nondegenerate
(symplectic), and the Poisson structure $\Pi_{trans}$ (the transverse
structure of $\Pi$ at 0) vanishes at 0. If $\Pi_{trans}$ is trivial, i.e.
the Poisson structure $\Pi$ is regular near $0$, then the problem local
normal forms of Hamiltonian vector fields near $0$ is reduced to the usual
problem of normal forms Hamiltonian vector fields (with parameters) on a
symplectic manifold. Here we are interested in the case when $\Pi_{trans}$
is not trivial. We will restrict our attention to the case when the linear
part of $\Pi_{trans}$ corresponds to a semisimple Lie algebra $\fg$.
According to linearization theorems of Weinstein
\cite{Weinstein-Poisson1983} and Conn \cite{Conn-Analytic1984}, we may
identify $((\K^m,0),\Pi_{trans})$ with a neighborhood of 0 of the dual
$\fg^*$ of $\fg$ equipped with the associated linear (Lie-Poisson)
structure. In other words, there is a local system of coordinates
$(x_1,y_1,\hdots,x_l,y_l,z_1,\hdots,z_m)$ on $\K^{2l+m}$ such that
$\Pi_{symp} = \sum_{i=1}^l {\partial \over \partial x_i} \wedge {\partial
\over \partial y_i}$, $\Pi_{trans} = \Pi_\fg = {1 \over 2}\sum_{i,j,k}
c_{ij}^k z_k {\partial \over \partial z_i} \wedge {\partial \over \partial
z_j}$ with $c_{ij}^k$ being structural constants of $\fg$, and
\begin{equation}
\Pi = \sum_{i=1}^l {\partial \over \partial x_i} \wedge {\partial \over
\partial y_i} +  {1 \over 2}\sum_{i,j,k} c_{ij}^k z_k {\partial
\over \partial z_i} \wedge {\partial \over \partial z_j} .
\end{equation}
Such a coordinate system will be called a \emph{canonical coordinate
system} of $\Pi$ near $0$. Let $H$ be a formal or analytic function on
$((\K^n,0),\Pi)$. We will assume that the Hamiltonian vector field $X_H$
of $H$ vanishes at $0$. Note that the differential of $H$ does not
necessarily vanish at $0$ (for example, if $n =0$ then we always have
$X_H(0) = 0$ for any $H$). We may assume that $H(0) = 0$.

We have the following generalization of Birkhoff normal form
\cite{Birkhoff} (see Theorem \ref{thm:formal-Birkhoff}):

{\bf Theorem C}. {\it With the above notations and assumptions, there is a
formal canonical coordinate system $(\hat{x}_i,\hat{y}_i,\hat{z}_j)$, in
which $H$ satisfies the following equation:
$$
\{H,H_{ss}\}=0\ ,
$$
where $H_{ss}$ is a (nonhomogeneous quadratic) function such that its
Hamiltonian vector field $X_{H_{ss}}$ is linear and is the semisimple part
of the linear part $X^{(1)}_H$ of $X_H$ (in this coordinate system). In
particular, the semisimple part of the linear part of $X_H$ is a
Hamiltonian vector field.}

Note that the normalizing canonical coordinates given in the above theorem
are only formal in general. The problem of existence of a local analytic
normalization for a Hamiltonian vector field (even in the symplectic case)
is much more delicate than for a general vector field, due to
``auto-resonances'' (e.g, if $\lambda$ is an eigenvalue of a Hamiltonian
vector field then $-\lambda$ also is). However, there is one particular
situation where one knows that a local analytic normalization always
exists, namely when the Hamiltonian vector field is analytically
integrable. See \cite{Zung-Birkhoff} for the case of integrable
Hamiltonian vector fields on symplectic manifolds. Here we can generalize
the main result of \cite{Zung-Birkhoff} to our situation (see Theorem
\ref{thm:integrability-cvgtPB}):

{\bf Theorem D}. {\it Assume that $\K =\C$, the Hamiltonian function $H$
in Theorem C is locally analytic, and is analytically integrable in the
generalized Liouville sense. Then the normalizing canonical coordinate
system $(\hat{x}_i,\hat{y}_i,\hat{z}_j)$ can be chosen locally analytic.}

We conjecture that the above theorem remains true in the real case ($\K =
\R$). Recall (see, e.g., \cite{Zung-survey} and references therein) that a
Hamiltonian vector field $X_H$ on a Poisson manifold $(M,\Pi)$ of
dimension $n$ is called \emph{integrable in generalized Liouville sense}
if there are nonnegative integers $p,q$ with $p+q = n$, $p$ pairwise
commuting Hamiltonian functions $H_1,\hdots, H_p$ ($\{H_i,H_j\} = 0\
\forall i,j$) with $H_1 = H$ and $q$ first integrals $F_1,\hdots,F_p$,
such that $X_{H_i}(F_j) = 0 \ \forall\ i,j$, and $dF_1 \wedge \hdots
\wedge dF_q \neq 0$ and $X_{H_1} \wedge \hdots \wedge X_{H_p} \neq 0$
almost everywhere. (The Liouville case corresponds to $p = q = n/2$ and
$F_i = H_i$). \emph{Analytic} integrability means that all Hamiltonian
functions and vector fields in question are analytic.

\section{Homogeneous Poisson structures}
\label{section:Poisson-homo}

Following \cite{Marle-Jacobi}, we will use the following terminology: a
{\it homogeneous Poisson structure} on a manifold $M$ is a couple
$(\Pi,X)$ where $\Pi$ is a Poisson structure and $X$ a vector field which
satisfies the relation
\begin{equation}
[X,\Pi]=-\Pi \,, \label{def:homogeneous}
\end{equation}
where the bracket is the Schouten bracket.

\begin{rem} Poisson structures which satisfy the above condition are also called
{\it exact}, in the sense that the Poisson tensor is a coboundary in the
associated Lichnerowicz complex which defines Poisson cohomology. They
have nothing to do with another kind of homogeneous spaces, namely those
which admit a transitive group action.
\end{rem}

An analog of Weinstein's splitting theorem for homogeneous Poisson
structures is given in \cite{Marle-Jacobi}, and it reduces the study of
normal forms of homogeneous Poisson structures to the case when both $\Pi$
and $X$ vanish at a point. So we will assume that $(\Pi,X)$ is a
homogeneous Poisson structure defined in a neighborhood of 0 in $\K^n$,
where $\K = \R$ or $\C$, such that
\begin{equation}
\Pi(0)=0\quad {\mbox { and }}\quad X(0)=0\,.
\label{eqn:assumptionPoissHomo}
\end{equation}

We are interested in the linearization of these structures, i.e.
simultaneous linearization of $\Pi$ and $X$. Denote by $\Pi_1$ and
$X^{(1)}$ the linear parts of $\Pi$ and $X$ respectively. Then the terms
of degree 1 of Equation (\ref{def:homogeneous}) imply that
$(\Pi_1,X^{(1)})$ is again a homogeneous Poisson structure.

In this paper, we will assume that the linear Poisson structure $\Pi_1$
corresponds to a semisimple Lie algebra, which we denote by $\fg$. Then,
according to linearization results of Weinstein
\cite{Weinstein-Poisson1983} (for the formal case) and Conn
\cite{Conn-Analytic1984} (for the analytic case), the Poisson structure
$\Pi$ can be linearized. In other words, there is a local coordinate
system $(x_1,\hdots,x_n)$ on $(\K^n,0)$, in which
\begin{eqnarray}
\Pi = \Pi_1 &=& {1 \over 2}\sum_{ijk} c_{ij}^kx_k \frac{\partial
}{\partial
x_i}\wedge\frac{\partial }{\partial x_j}\, , {\mbox { or }}  \\
\{x_i,x_j\} &=& \sum_k c_{ij}^k x_k ,
\end{eqnarray}
where $c_{ij}^k$ are structural constants of $\fg$. In order to linearize
$(\Pi,X) = (\Pi_1,X)$, it remains to linearize $X$ by local (formal or
analytic) diffeomorphisms which preserve the linear Poisson structure
$\Pi_1$.

\subsection{Formal linearization} First consider the complex case ($\K = \C$).
Let $X$ be a formal vector field on $\C^n$ such that $(\Pi_1,X)$
forms a homogeneous Poisson structure on $\C^n$.
Denote by
\begin{equation}
I=\sum_{i=1}^n x_i \frac{\partial }{\partial x_i}
\end{equation}
the Euler vector field written in coordinates $(x_1,\hdots,x_n)$. Since
this vector field satisfies the relation $[I,\Pi_1]=-\Pi_1$, we can write
$X$ as
\begin{equation}
X= I+Y\,,
\end{equation}
where $Y$ is a Poisson vector field with respect to $\Pi_1$, i.e.,
$[Y,\Pi_1]=0$. It is well-known that, since the complex Lie algebra $\fg$
is semisimple by assumptions, the first formal Poisson cohomology space of
$\Pi_1$ is trivial (see, e.g., \cite{DZ-PoissonBook}), i.e. any formal
Poisson vector field is Hamiltonian. In particular, we have
\begin{equation}
Y = X_h = - [h, \Pi_1]
\end{equation}
for some formal function $h$. Writing the Taylor expansion
$h=h_1+h_2+h_3+\cdots$ where each $h_r$ is a polynomial of degree $r$, we
have
\begin{equation}
X= I+X_{h_1}+X_{h_2}+X_{h_3}+\cdots
\end{equation}

Denote by $X^{(1)}= I+X_{h_1}$ the linear part of $X$. In order to
linearize $X$ (while preserving the linearity of $\Pi = \Pi_1$), we want
to kill all the terms $X_{h_r}$ with $r\geq 2$, using a sequence of
changes of coordinates defined by flows of Hamiltonian vector fields with
respect to $\Pi_1$. Working degree by degree, we want to find for each $r$
a homogeneous polynomial $g_r$ of degree $r$ such that
\begin{equation} \label{eqn:homotopy1}
[X^{(1)},X_{g_r}]=X_{h_r}\,.
\end{equation}

Note that $[X_{h_1},X_{g_r}] = X_{\{h_1,g_r\}}$, and $[I, X_{g_r}] = (r-1)
X_{g_r}$ because $X_{g_r}$ is homogeneous of degree $r$. Hence Relation
(\ref{eqn:homotopy1}) will be satisfied if $g_r$ satisfies the following
relation:
\begin{equation}
(r-1) g_r+\{h_1,g_r\}=h_r\,.
\end{equation}

Remark that, $h_1$ can be viewed as an element of $\fg$, and $h_r$, $g_r$
may be identified with elements of the symmetric power $S^r(\fg)$ of
$\fg$. Under this identification, $\{h_1,g_r\}$ is nothing but the result
of the adjoint action of $h_1 \in \fg$ on $g_r \in S^r(\fg)$.

We will suppose that $h_1$ is a semisimple element of $\fg$, and denote by
$\fh$ a Cartan subalgebra of $\fg$ which contains $h_1$. According to the
root decomposition of $\fg$ with respect to $\fh$, we can choose a basis
$(x_1,\hdots,x_n)$ of $\fg$, and elements $\alpha_1,\hdots,\alpha_n$ of
$\fh^*$, such that
\begin{equation}
[y,x_i] = \langle \alpha_i, y\rangle x_i \ \ \forall\ y \in \fh,\ \forall
i=1,\hdots,n.
\end{equation}
Each $\alpha_i$ is either 0 (in which case $x_i \in \fh$) or a root of
$\fg$ (in which case $x_i$ belongs to the root subspace $\fg_{\alpha_i}$
of $\fg$).


We define for each $r\geq 2$ the linear operator
\begin{eqnarray*}
\Theta_r : S^r(\fg) & \longrightarrow & S^r(\fg)\\
a & \longmapsto & (r-1) a+\{ h_1,a \}\ .
\end{eqnarray*}
Each monomial $\prod_i x_i^{\lambda_i}$ of degree
$|\lambda|=\sum\lambda_i=r$ is an eigenvector of this linear operator:
\begin{equation}
\Theta_r(\prod_i x_i^{\lambda_i})= \big(r - 1 + \sum_{i=1}^n \lambda_i
\langle \alpha_i, h_1 \rangle) \prod_i x_i^{\lambda_i}\,.
\label{eqn:thetar}
\end{equation}

\begin{defn} With the above notations, we will say that
$(\Pi_1,X^{(1)})$ is {\it semisimple nonresonant} if $h_1$ is a semisimple
element of $\fg$ and the eigenvalues of $\Theta_r$ don't vanish, i.e., for
any $r \geq 2$ and any $(\lambda_1,\hdots,\lambda_n) \in \Z_+^n$ such that
$\sum \lambda_i = r$ we have $r - 1 + \sum_{i=1}^n \lambda_i \langle
\alpha_i, h_1 \rangle \neq 0$.
\end{defn}

\begin{rem}
It is easy to see that the above nonresonance condition is a generic
position condition, and the subset of elements which do not satisfy this
condition is of codimension 1. In fact, if the Cartan subalgebra $\fh$ is
fixed, then the set of elements $h_1 \in \fh$ such that $(\Pi_1, I +
X_{h_1})$ is resonant is a countable union of affine hyperplanes in $\fh$
which do not contain the origin, and there is a neighborhood of $0$ in
$\fh$ such that if $h_1$ belongs to this neighborhood then
$(\Pi_1,X^{(1)})$ is automatically nonresonant.
\end{rem}

{\it The algorithm of formal linearization}. We now show how to linearize
$(\Pi_1,X)$, by killing the nonlinear terms of $h$ step by step, provided
that $(\Pi_1,X^{(1)})$ is nonresonant. Actually, at each step, we will
kill not just one term $h_d$, but a whole block of  $2^d$ consecutive
terms. This ``block killing'' will be important in the next section when
we want to show that, under some Diophantine-type condition, our formal
linearization process actually yields a local analytic linearization.

For each $q\geq 0$, denote by $\hat{\mathcal{O}}_q$ the space of formal
power series on $\C^n$ of order greater or equal to $q$, i.e. without
terms of degree $< q$.

We begin with $X=X^{(1)}\; {\mbox { mod }}\,\hat{\mathcal{O}}_2$, and will
construct a sequence of formal vector fields ${(X_d)}_d$ and
diffeomorphisms ${(\varphi_d)}_d$, such that $X_0 =X$ and, for all $d\geq
0$,
\begin{eqnarray}
X_d &=& X^{(1)}\;{\mbox { mod }}\,\hat{\mathcal{O}}_{2^{d}+1}\ , \\
X_{d+1} &=& {\varphi_d}_\ast X_d.
\end{eqnarray}
Assuming that we already have $X_d$ for some $d \geq 0$, we will construct
$\varphi_d$ (and $X_{d+1} = {\varphi_d}_\ast X_d$). We write
\begin{equation}
X_d=X^{(1)}+X_{H_d} \; {\mbox { mod }}\,\hat{\mathcal{O}}_{2^{d+1}+1}\,,
\label{eqn:recurenced}
\end{equation}
where $H_d$ is a polynomial of degree $\leq 2^{d+1}$ in
$\hat{\mathcal{O}}_{2^d+1}$, i.e. $H_d$ is a sum of homogeneous
polynomials of degrees between $2^d+1$ and $2^{d+1}$ (we also could write
abusively that $H_d$ is in
$\hat{\mathcal{O}}_{2^d+1}/\hat{\mathcal{O}}_{2^{d+1}+1}$). Under the
nonresonance condition, there exists a polynomial $G_d$ of order $\geq
2^d+1$ and of degree $\leq 2^{d+1}$ such that if we write $G_d=
G_d^{(2^d+1)}+\hdots+ G^{(2^{d+1})}_d$ (where $G^{(u)}_d$ is homogeneous
of degree $u$) and the same for $H_d$, we have
\begin{equation} \label{eqn:Gdu}
\Theta_u(G^{(u)}_d)=(u-1) G^{(u)}_d+ \{h_1,G^{(u)}_d \}=H^{(u)}_d \,,
\end{equation}
for every $u\in\{2^d+1,\hdots,2^{d+1}\}$. It implies that we have
\begin{equation}
[X^{(1)},X_{G_d}]=X_{H_d}\, , \label{eqn:primitived}
\end{equation}
where $X_{G_d}$ denotes the Hamiltonian vector field of $G_d$ with respect
to $\Pi_1$ as usual. Now, we define the diffeomorphism $\varphi_d=\exp
X_{G_d}$ to be the time-1 flow of $X_{G_d}$. We then have
\begin{equation}
X_{d+1}:={\varphi_d}_\ast X= X^{(1)} + X_{H_{d+1}}\; {\mbox { mod }}\,
\hat{\mathcal{O}}_{2^{d+2}+1}\,,
\end{equation}
where $H_{d+1}$ is a polynomial of degree $2^{d+2}$ in
$\hat{\mathcal{O}}_{2^{d+1}+1}$.

Constructed in this way, it is clear that the successive compositions of
the diffeomorphisms $\varphi_d$ converge in the formal category to a
formal diffeomorphism $\Phi_\infty$ which satisfies ${\Phi_\infty}_\ast
X=X^{(1)}$ and which preserves the linear Poisson structure $\Pi_1$.

Consider now the real case ($\K = \R$, and $\fg$ is a real semisimple Lie
algebra). By complexication, we can view real objects as holomorphic
objects with real coefficients, and then repeat the above algorithm. In
particular, under the nonresonance condition, we will find homogeneous
polynomials $G_d^{(u)}$ which satisfy Equation (\ref{eqn:Gdu}), i.e.,
\begin{equation}
(u-1) G^{(u)}_d+ \{h_1,G^{(u)}_d \}=H^{(u)}_d.
\end{equation}
Remark that, in the real case, the operator $\Theta_u: G^{(u)} \mapsto
(u-1) G^{(u)} + \{h_1,G^{(u)} \}$ is real (and is invertible under the
nonresonance condition), and $H^{(u)}_d$ is real, so $G^{(u)}_d$ is also
real. This means that the coordinate transformations constructed above are
real in the real case.

We have proved the following:
\begin{thm}
Let $(\Pi,X)$ be a formal homogeneous Poisson structure on $\K^n$ (where
$\K$ is $\C$ or $\R$) such that the linear part $\Pi_1$ of $\Pi$
corresponds to a semisimple Lie algebra. Assume that its linear part
$(\Pi_1,X^{(1)})$ is semisimple nonresonant. Then there exists a formal
diffeomorphism which sends $(\Pi,X)$ to $(\Pi_1,X^{(1)})$.
\label{thm:Poisson-hom-formal}
\end{thm}

\subsection{Analytic linearization} Now we work in the local
analytic context, i.e. the vector field $X$ is supposed to be analytic on
$(\K^n,0)$. In order to show that the algorithm given in the previous
subsection leads to a local analytic linearization, in addition to the
nonresonance condition we will need a Diophantine-type condition, similar
to Bruno's $\omega$-condition for the analytic linearization of vector
fields \cite{Bruno,Bruno2}.

Keeping the notations of the previous subsection, for each $d\geq 1$, put
\begin{equation}\label{eqn:omega_d}
\omega_d=\min \left\{ {1 \over 2d} , \min\Big\{\Big| |\lambda| - 1
+ \sum_{i=1}^n \lambda_i\langle \alpha_i,h_1\rangle \Big| \;;\;
\lambda\in\Z_+^n \,{\mbox { and }}\, 2 \leq |\lambda|\leq 2^{d+1}
\Big\} \right\}\, .
\end{equation}

\begin{defn}\label{defn:omega_condition}
 We will say that $X^{(1)}$, or more precisely that a semisimple
nonresonant linear homogeneous
Poisson structure $(\Pi_1,X^{(1)})$ satisfies the
\emph{$\omega$-condition} if
\begin{equation}\label{eqn:Bruno}
\sum_{d=1}^\infty \frac{-\log\omega_d}{2^d}<\infty\,.
\end{equation}
\end{defn}

Remark that, similarly to other situations, the set of $X^{(1)}$ which
satisfy the about $\omega$-condition is of full measure. More precisely,
we have:

\begin{prop} \label{prop:fullmeasure}
The set of elements $h$ of a given Cartan subalgebra $\fh$ such that
$X^{(1)} = I + X_h$ does not satisfy the $\omega$-condition
(\ref{eqn:Bruno}) is of measure 0 in $\fh$.
\end{prop}

See the Appendix for a straightforward proof of the above proposition.

Using the same analytical tools as in the proof of Bruno's theorems about
linearization of analytic vector fields \cite{Bruno,Bruno2}, we will show
the following theorem:
\begin{thm}
Let $(\Pi,X)$ be an analytic homogeneous Poisson structure on $(\K^n,0)$
(where $\K$ is $\C$ or $\R$) such that the linear part $\Pi_1$ of $\Pi$
corresponds to a semisimple Lie algebra. Suppose that its linear part
$(\Pi_1,X^{(1)})$ is semisimple nonresonant and satisfies the
$\omega$-condition. Then there exists a local analytic diffeomorphism
which sends $(\Pi,X)$ to $(\Pi_1,X^{(1)})$.
\label{thm:Poisson-hom-analytic}
\end{thm}

\emph{Proof}. Due to Conn's theorem \cite{Conn-Analytic1984}, we can
assume that $\Pi = \Pi_{1}$ is already linear. The process to linearize
the vector field $X$ is the same as in the formal case, noting that if we
start with an analytic vector field, the diffeomorphisms $\varphi_d$ that
we constructed will be analytic too (as is the vector fields $X_d$). We
just have to check the convergence of the sequence
$\Phi_d=\varphi_d\circ\hdots\circ\varphi_1$ in the analytic setup.

We will assume that $\K = \C$ (the real case can be reduced to the complex
case by the same argument as given in the previous subsection). Denote by
$\mathcal{O}_q$  the vector space of local analytic functions of
$(\K^n,0)$ of order greater or equal to $q$ (i.e. without terms of degree
$<q$).

For each positive real number $\rho>0$, denote by $D_\rho$ the ball
$\{x=(x_1,\hdots,x_n)\in\C^n\,;\, |x_i|<\rho\}$ and if
$f=\sum_{\lambda\in\N^n} a_\lambda x^\lambda$ is an analytic function on
$D_\rho$ we define the following norms:
\begin{eqnarray}
|f|_\rho &:=& \sum_\lambda |a_\lambda| \rho^{|\lambda|}\ , \\
\|f\|_\rho &:=& \sup_{z\in D_\rho} |f(z)| \ .
\end{eqnarray}
In the same way, if $F=(F_1,\hdots,F_n)$ is a vector-valued local map then
we put $|F|_\rho := \max \{|F_1|_\rho,\hdots,|F_n|_\rho\}$ and similarly
for $||F||_\rho$. These norms satisfy the following properties.

\begin{lem}\label{lem:normes}
Let $\rho$ and $\rho^\prime$ be two real numbers such that
$0<\rho^\prime<\rho$. If $f\in \mathcal{O}_q$ is an analytic
function on $D_\rho$, then \\
a)
\begin{equation}
\|f\|_\rho \leq |f|_\rho \quad {\mbox { and }} \quad
|f|_{\rho^\prime}\leq \frac{1}{1-(\rho^\prime /\rho)}\|f\|_\rho\,.
\label{eqn:eqnormes}
\end{equation}
b)
\begin{equation}
|f|_{\rho^\prime}<\big(\frac{\rho^\prime}{\rho}\big)^q |f|_\rho\,.
\label{eqn:chgtrad}
\end{equation}
c) Let $R > 0$ be a positive constant. Then there is a natural
number $N$ such that for any $d > N$, if $q = 2^d+1$,
${(\frac{1}{(2d)(2^d)})}^{1/(2^d+1)}\rho = \rho' \geq R$, and
$f\in \mathcal{O}_q$ is an analytic function on $D_\rho$, then we
have
\begin{equation} \label{eqn:normederiv}
|df|_{\rho^\prime} \leq |f|_\rho \,.
\end{equation}
\end{lem}

The proof of the above lemma is elementary (see the Appendix).

It is important to remark that, with the same notations as in the formal
case, for $\rho>0$, we have, by (\ref{eqn:thetar}):
\begin{equation} \label{eqn:smalldiv}
|X_{G_d}|_\rho \leq \frac{1}{\omega_d} |X_{H_d}|_\rho\,.
\end{equation}

Put $\rho_0 =1$, and define the following two decreasing sequences of
radii ${(r_d)}_d$ and ${(\rho_d)}_d$ by
\begin{eqnarray}
r_d &:=& \big(\frac{\omega_d}{2^d}\big)^{1/(2^d+1)} \rho_{d-1}\ ,\\
\rho_d &:=& (1-\frac{1}{d^2}) r_d\ .
\end{eqnarray}
We have
\begin{equation}
\hdots<\rho_{d+1}<r_{d+1}<\rho_d<r_d<\rho_{d-1}<\hdots ,
\end{equation}
and it is clear, by the $\omega$-condition (\ref{eqn:Bruno}), that the
sequences ${(r_d)}_d$ and ${(\rho_d)}_d$ converge to a strictly positive
limit $R>0$. Moreover, they satisfy the following properties:

\begin{lem}
For $d$ sufficiently large, we have
\begin{itemize}
\item[a)] $r_d-\rho_d > \frac{1}{2^d}$ , \\
\item[b)] $\rho_d-r_{d+1}>\frac{1}{2^d}$ .
\end{itemize}
\label{lem:radii}
\end{lem}

The proof of Lemma \ref{lem:radii} is elementary (see the Appendix).

\begin{lem}
For $d$ sufficiently large, if $|X_d-X^{(1)}|_{\rho_{d-1}}<1$,
then
\begin{equation}
D_{r_{d+1}} \subset \varphi_d (D_{\rho_d}) \subset D_{r_d}\,,
\end{equation}
and moreover, we have $|{\varphi_d}_\ast
X_d-X^{(1)}|_{\rho_{d}}<1$. \label{lem:boules}
\end{lem}

{\it Proof }. $\bullet$ We first prove the second inclusion : $\varphi_d
(D_{\rho_d}) \subset D_{r_d}$.

We have
\begin{equation}
X_d=X^{(1)} + X_{H_d} \; {\mbox { mod }} \mathcal{O}_{2^{d+1}+1}
\end{equation}
where $H_d$ is a polynomial formed by homogenous terms of degree
between $2^d+1$ and $2^{d+1}$. By (\ref{eqn:smalldiv}), we write
\begin{equation}
|X_{G_d}|_{\rho_{d-1}}<\frac{1}{\omega_d}|X_{H_d}|_{\rho_{d-1}}\,.
\end{equation}
Then, by (\ref{eqn:chgtrad}), we get
\begin{equation}
|X_{G_d}|_{r_d}<\frac{1}{\omega_d}\big(\frac{r_d}{\rho_{d-1}}\big)^{2^d+1}
|X_{H_d}|_{\rho_{d-1}}
\end{equation}
And, using the assumption $|X_d-X^{(1)}|_{\rho_{d-1}}<1$, we
obtain
\begin{equation}
|X_{G_d}|_{r_d}<\frac{1}{2^d}\,. \label{eqn:XGrd}
\end{equation}
Finally, Lemma \ref{lem:radii} gives
\begin{equation}
\|X_{G_d}\|_{r_d}\leq |X_{G_d}|_{r_d}<r_d-\rho_d\,,
\label{eqn:Xrrho}
\end{equation}
which implies the inclusion $\varphi_d (D_{\rho_d}) \subset
D_{r_d}$.\\

$\bullet$ Now, we prove the first inclusion $D_{r_{d+1}} \subset \varphi_d
(D_{\rho_d})$. For any $x$ on the boundary $S_{\rho_d}$ of $D_{\rho_d}$,
we define $x_1:=\frac{r_d+\rho_d}{2}\frac{x}{|x|}$ and
$x_2:=r_d\frac{x}{|x|}$. We construct a map $\hat\phi_d :
D_{r_d}\longrightarrow D_{r_d}$ which is $\varphi_d$ on $D_{\rho_d}$ and
defined on $D_{r_d}\backslash D_{\rho_d}$ by the following : for $\mu\in
[0,1]$ and $x\in S_{\rho_d}$ we put
\begin{eqnarray*}
\hat\phi_d(\mu x+(1-\mu)x_1) &=& \mu\varphi_d(x) + (1-\mu)x\\
\hat\phi_d(\mu x_1+(1-\mu)x_2) &=& \mu x + (1-\mu)x_2
\end{eqnarray*}
This map is continuous and is the identity on the boundary of $D_{r_d}$
thus, by Brouwer's theorem, $\hat\phi_d(D_{r_d})=D_{r_d}$.

Let $x$ be an element of the boundary $S_{\rho_d}$ of
$D_{\rho_d}$.

If $z=\mu x +(1-\mu) x_1$ (for $\mu\in [0,1]$) then we have
\begin{eqnarray*}
|\hat\phi_d(z)| &=& |x+\mu (\varphi_d(x)-x)|\\
           &\geq& |x|-\mu |\varphi_d(x)-x|
\end{eqnarray*}
Now, we write $|\varphi_d(x)-x|\leq \big|\int_0^1
X_{G_d}\big(\varphi_d^t(x)\big)dt\big|$ where $\varphi_d^t$ is the
flow of $X_{G_d}$. As above, according to (\ref{eqn:Xrrho}),
$\varphi_d^t(x)$ is in $D_{r_d}$ for all $t\in[0,1]$ and then, we
get $|\varphi_d(x)-x|\leq \|X_{G_d}\|_{r_d}<\frac{1}{2^d}$.
Therefore, by Lemma \ref{lem:radii}, we get
\begin{equation}
|\hat\phi_d(z)|>r_{d+1}\,.
\end{equation}

Now, if $z=\mu x_1 +(1-\mu) x_2$ ($\mu\in [0,1]$) then we have
\begin{equation}
|\hat\phi_d(z)|=\big(\mu+(1-\mu)\frac{r_d}{|x|}\big) |x|\geq
|x|>r_{d+1}\,.
\end{equation}
As a conclusion, if $y$ is in $D_{r_{d+1}}$ then, by the surjectivity of
$\hat\phi_d$, $y=\hat\phi_d(z)$ with, a priori, $z$ in $D_{r_d}$. We saw
above that in fact $z$ cannot be in $D_{r_d}/D_{\rho_d}$. Therefore, since
$\hat\phi_d=\varphi_d$ on $D_{\rho_d}$, we get $y=\varphi_d(z)$ with $z$
in
$D_{\rho_d}$.\\

$\bullet$ Finally, we check that $|{\varphi_d}_\ast X_d -
X^{(1)}|_{\rho_d}<1$. We write the obvious inequality
\begin{equation}
|{\varphi_d}_\ast X_d - X^{(1)}|_{\rho_d} \leq |{\varphi_d}_\ast
X_d - X_d|_{\rho_d} + |X_d-X^{(1)}|_{\rho_d}\,.
\end{equation}
By (\ref{eqn:chgtrad}), we have
\begin{equation}
|X_d-X^{(1)}|_{\rho_d} <
\big(\frac{\rho_d}{\rho_{d-1}}\big)^{2^d+1} <
\frac{\omega_d}{2^d}{\big(1-\frac{1}{d^2}\big)}^{2^d+1}\,.
\end{equation}

Now, we just have to estimate the term $|{\varphi_d}_\ast X_d -
X_d|_{\rho_d}$.  To do that, we use the inequalities of Lemma \ref{lem:normes}. 
The drawback of these inequalities is that they sometimes induce a change of radius. 
Therefore, we define the following intermediar radii (between $\rho_d$ and $r_d$) :
\begin{eqnarray*}
\rho_d^{(1)} &=& \rho_d (1+\frac{1}{5d^2}) \\
\rho_d^{(2)} &=& \rho_d^{(1)} + \frac{3}{2}\frac{1}{2^d} \\
\rho_d^{(3)} &=& \rho_d^{(2)}\big( (2d)(2^d) \big)^{\frac{1}{2^d+1}} \\
\rho_d^{(4)} &=& \rho_d^{(3)} (1+\frac{1}{5d^2})
\end{eqnarray*}

Let us explain a little bit the definitions of these radii :
\begin{itemize}
\item[- ] $\rho_d^{(1)}$ (resp. $\rho_d^{(4)}$) is defined from $\rho_d$ (resp. $\rho_d^{(3)}$) 
in order to use inequality (\ref{eqn:eqnormes}) and have 
$$
\frac{1}{1-\frac{\rho_d}{\rho_d^{(1)}}}\sim 5d^2
$$
which does not grow too quickly.
\item[- ] $\rho_d^{(2)}$ is defined in order to have (recall  (\ref{eqn:XGrd}))
\begin{equation}
\rho_d^{(2)}-\rho_d^{(1)}>\frac{1}{2^d}>\parallel X_{G_d} \parallel_{r_d}
\label{eqn:a}
\end{equation}
\item[- ] $\rho_d^{(3)}$ is defined in order  to use inequality (\ref{eqn:normederiv}).
\item[- ] Finally, {\it if $d$ is sufficiently large}, the differences $\rho_d^{(1)}-\rho_d$, 
$\rho_d^{(2)}-\rho_d^{(1)}$, $\rho_d^{(3)}-\rho_d^{(2)}$ and $\rho_d^{(4)}-\rho_d^{(3)}$
are strictly smaller than $\frac{r_d}{5d^2}$ and then, 
\begin{equation}
r_d-\rho_d^{(4)}>\frac{r_d}{d^2}-\frac{4r_d}{5d^2}>
\frac{r_d}{5d^2}>\frac{1}{2^d}>\parallel X_{G_d} \parallel_{r_d}\,.
\label{eqn:b}
\end{equation}
\end{itemize}

We have, by (\ref{eqn:eqnormes}),
\begin{equation}
|{\varphi_d}_\ast X_d - X_d|_{\rho_d}\leq
\frac{1}{1-\frac{\rho_d}{\rho_d^{(1)}}} \|{\varphi_d}_\ast X_d -
X_d\|_{\rho_d^{(1)}} = (5d^2+1)\|{\varphi_d}_\ast X_d -
X_d\|_{\rho_d^{(1)}}\,. \label{est1}
\end{equation}
If $x$ is in $D_{\rho_d^{(1)}}$ then we have
\begin{eqnarray}
\big|\big({\varphi_d}_\ast X_d - X_d\big)(x)\big| &=& 
    \Big|\int_0^1 {\varphi_d^t}_\ast [X_{G_d},X_d] (x) dt\Big|\\
 &=& \Big|\int_0^1 \Big( d\varphi_d^t([X_{G_d},X_d])\Big) \big(\varphi_d^{-t}(x)\big) dt \Big|
\,.\nonumber
\label{est2}
\end{eqnarray}
Since $\parallel X_{G_d}\parallel_{\rho_d^{(2)}} \leq \parallel X_{G_d}\parallel_{r_d} <
\rho_d^{(2)}-\rho_d^{(1)}$  (by (\ref{eqn:a})), 
$\varphi_d^{-t}(x)$ belongs to $D_{\rho_d^{(2)}}$ for all $t\in [0,1]$.
We  then get 
\begin{equation}
\|{\varphi_d}_\ast X_d - X_d\|_{\rho_d^{(1)}} \leq 
\int_0^1 \| d\varphi_d^t([X_{G_d},X_d]) \|_{\rho_d^{(2)}} dt\,.
\label{eqn:c}
\end{equation}
We can write $\varphi_d^t=Id+\xi_d^t$ where the $n$ components of $\xi_d^t$ are 
functions in $\mathcal{O}_{2^d+1}$. We have the estimates
\begin{eqnarray*}
\|d\xi_d^t\|_{\rho_d^{(2)}} & \leq & | d\xi_d^t |_{\rho_d^{(2)}} \quad {\mbox { by (\ref{eqn:eqnormes})}}\\
& \leq & | \xi_d^t |_{\rho_d^{(3)}} \quad {\mbox { by (\ref{eqn:normederiv}) }}\\
& \leq & (5d^2+1) \| \xi_d^t \|_{\rho_d^{(4)}} \quad {\mbox { by (\ref{eqn:eqnormes})}}\,.
\end{eqnarray*}

If $x$ is in $D_{\rho_d^{(4)}}$ then we can write
\begin{equation}
\xi_d^t(x)=\varphi_d^t(x)-x=\int_0^t X_{G_d}(\varphi_d^u (x)) du\,.
\end{equation}
Since $\|X_{G_d}\|_{r_d}<\frac{1}{2^d}<r_d-\rho_d^{(4)}$ (see (\ref{eqn:b})), we have
$\varphi_d^u (x) \in D_{r_d}$ for all $u$ in $[0,t]$. 
Thus $\|\xi_d^t\|_{\rho_d^{(4)}} \leq \|X_{G_d}\|_{r_d}<\frac{1}{2^d}$ which gives
\begin{equation}
\|d\xi_d^t\|_{\rho_d^{(2)}} \leq \frac{5d^2+1}{2^d}\,,
\end{equation}
and then, by (\ref{eqn:c}),
\begin{equation}
\|{\varphi_d}_\ast X_d - X_d\|_{\rho_d^{(1)}} \leq
\big( 1+\frac{5d^2+1}{2^d}\big) \|[X_{G_d},X_d]) \|_{\rho_d^{(2)}}\,.
\end{equation}
We then deduce by (\ref{est1}) that 
\begin{equation}
|{\varphi_d}_\ast X_d - X_d|_{\rho_d}\leq (5d^2+1) 
\big( 1+\frac{5d^2+1}{2^d}\big) \|[X_{G_d},X_d]) \|_{r_d}\,.
\end{equation}

Finally, we just have to estimate $\|[X_{G_d},X_d]) \|_{r_d}$. We first have by (\ref{eqn:eqnormes}),
\begin{equation}
\|[X_{G_d},X_d]) \|_{r_d}\leq |[X_{G_d},X_d]) |_{r_d}\,.
\end{equation}

Now, we write
\begin{equation}
[X_{G_d},X_d]=[X_{G_d},X^{(1)}]+[X_{G_d},X_d-X^{(1)}]
=-X_{H_d}+[X_{G_d},X_d-X^{(1)}]\,,
\end{equation}
which gives, by (\ref{eqn:normederiv}), recalling that
$\omega_d\leq\frac{1}{2d}$,
\begin{eqnarray*}
|{\varphi_d}_\ast X_d - X_d|_{\rho_d} &\leq&
(2d^2-1)\Big[|X_{H_d}|_{r_d}+|X_{G_d}|_{r_d}
|X_d-X^{(1)}|_{\rho_{d-1}}\\
   & & \quad + |X_{G_d}|_{\rho_{d-1}}
|X_d-X^{(1)}|_{r_d}\Big]\,.
\end{eqnarray*}
Using (\ref{eqn:smalldiv}) and (\ref{eqn:chgtrad}), we get
\begin{equation}
|X_{G_d}|_{r_d} |X_d-X^{(1)}|_{\rho_{d-1}} <
\frac{1}{\omega_d}\big(\frac{r_d}{\rho_{d-1}}\big)^{2^d+1}
|X_{H_d}|_{\rho_{d-1}} |X_d-X^{(1)}|_{\rho_{d-1}}\,,
\end{equation}
and then
\begin{equation}
|X_{G_d}|_{r_d} |X_d-X^{(1)}|_{\rho_{d-1}} < \frac{1}{2^d}\,.
\end{equation}
In the same way, one can prove that
\begin{equation}
|X_{G_d}|_{\rho_{d-1}} |X_d-X^{(1)}|_{r_d} < \frac{1}{2^d}\,.
\end{equation}
In addition, by (\ref{eqn:chgtrad}), we get
\begin{equation}
|X_{H_d}|_{r_d}\leq \big(\frac{r_d}{\rho_{d-1}}\big)^{2^d+1}
|X_{H_d}|_{\rho_{d-1}} <\frac{\omega_d}{2^d}\,.
\end{equation}
We deduce finally that
\begin{equation}
|{\varphi_d}_\ast X_d-X_d|_{\rho_d} <
(5d^2+1)\big( 1+\frac{5d^2+1}{2^d}\big) \big(\frac{\omega_d}{2^d}+\frac{2}{2^d}\big)\,.
\end{equation}
This gives the following estimate
\begin{equation}
|{\varphi_d}_\ast X_d - X^{(1)}|_{\rho_d} <
 (5d^2+1)\big( 1+\frac{5d^2+1}{2^d}\big) \big(\frac{\omega_d}{2^d}+\frac{2}{2^d}\big) +
\frac{\omega_d}{2^d}{\big(1-\frac{1}{2^d}\big)}^{2^d+1}\,,
\end{equation}
and the conclusion follows.\QED \\

{\it End of the proof of Theorem  \ref{thm:Poisson-hom-analytic}}. Let
$d_0$ be a positive integer such that Lemmas \ref{lem:radii} and
\ref{lem:boules} are satisfied for $d \geq d_0$. By the homothety trick
(dilate a given coordinate system by appropriate linear transformations),
we can assume that $|X_{d_0}-X^{(1)}|_{\rho_{d_0-1}}<1$.

By recurrence, for all $d\geq d_0$, we have
$$
D_{r_{d+1}} \subset \varphi_d (D_{\rho_d}) \subset D_{r_d}\,,
$$
which give
\begin{equation}
{\varphi_d}^{-1} (D_{r_{d+1}}) \subset D_{\rho_d}\,,
\label{eqn:invboules}
\end{equation}
for all $d\geq d_0$.

We consider the sequence ${(\Psi_d)}_d$ given by
$$
\Psi_d:=\varphi_0^{-1}\circ\varphi_1^{-1}\circ\hdots\circ\varphi_d^{-1}\,.
$$
Let $x$ be an element of $D_R$ ; recall that $R>0$ is the limit of
the decreasing sequences ${(r_d)}_d$ and ${(\rho_d)}_d$. Then $x$
belongs to the ball $D_{r_{d+1}}$ for any $d$ and if $d>d_0$, we
get by (\ref{eqn:invboules}), $\varphi_d^{-1}(x)\in
D_{\rho_d}\subset D_{r_d}$. In the same way, we get
$\varphi_{d-1}^{-1}\big(\varphi_d^{-1}(x)\big)\in
D_{\rho_{d-1}}\subset D_{r_{d-1}}$ and iterating this process, we
obtain
\begin{equation}
\varphi_{d_0}^{-1}\big(
\varphi_{d_0+1}^{-1}\circ\hdots\circ\varphi_d^{-1}(x)\big)\in
D_{r_{d_0}}\,.
\end{equation}
If we put $M=\sup_{z\in D_{d_0}}
\big|\varphi_{0}^{-1}\circ\hdots\circ\varphi_{d_0-1}^{-1}(z)\big|$,
we then obtain, for all $x$ in $D_R$ and all $d>d_0$,
\begin{equation}
\big|\Psi_d(x)\big|\leq M\,.
\end{equation}
The theorem follows.\QED

\section{Hamiltonian vector fields on Poisson manifolds}
\label{section:Hamiltonian}

In this section, we study normal forms of formal or analytic Hamiltonian
vector fields in the neighborhood of the origin on the Poisson manifold
$(\K^{2l+m},\Pi)$, where
\begin{equation} \label{Pi_canonical}
\Pi = \Pi_{symp} + \Pi_{trans} = \sum_{i=1}^l {\partial \over \partial
x_i} \wedge {\partial \over
\partial y_i} +  {1 \over 2}\sum_{i,j,k} c_{ij}^k z_k {\partial
\over \partial z_i} \wedge {\partial \over \partial z_j} .
\end{equation}
Here $\Pi_{symp} = \sum_{i=1}^l {\partial \over \partial x_i} \wedge
{\partial \over \partial y_i}$ is the standard symplectic Poisson
structure on $\K^{2n}$, and $\Pi_{trans} = \Pi_\fg = {1 \over
2}\sum_{i,j,k} c_{ij}^k z_k {\partial \over \partial z_i} \wedge {\partial
\over \partial z_j}$ is the associated linear Poisson structure on the
dual of a given semisimple Lie algebra $\fg$ of dimension $m$ over $\K$.

Let $H:(\K^{2l+m},0) \to (\K,0)$ be a formal or local analytic function
with $H(0)=0$, and consider the Hamiltonian vector field $X_H$ of $H$ with
respect to the above Poisson structure $\Pi = \Pi_{symp} + \Pi_\fg$. If
$X_H(0) \neq 0$, then it is well-known that it can be rectified, i.e.
there is a local canonical coordinate system
$(x_1,y_1,\hdots,x_l,y_l,z_1,\hdots,z_m)$ in which $H = x_1$ and $X_H =
{\partial \over \partial y_1}$. Here we will assume that $X_H(0) = 0$

\subsection{Formal Poincaré-Birkhoff normalization}
In this subsection, we will show that the vector field $X_H$ can be put
formally into Poincaré-Birkhoff normal form. More precisely, we have:

\begin{thm} \label{thm:formal-Birkhoff}
With the above notations, for any formal or local analytic function
$H:(\K^{2l+m},0) \to (\K,0)$, there is a formal canonical coordinate
system $(\hat{x}_i,\hat{y}_i,\hat{z}_j)$, in which the Poisson structure
$\Pi$ has the form
\begin{equation}
\Pi  = \sum_{i=1}^l {\partial \over \partial \hat{x}_i} \wedge {\partial
\over
\partial \hat{y}_i} +  {1 \over 2}\sum_{i,j,k} c_{ij}^k \hat{z}_k {\partial
\over \partial \hat{z}_i} \wedge {\partial \over \partial \hat{z}_j} ,
\end{equation}
and in which we have
\begin{equation} \label{eqn:defPB}
H = H_{ss} + \tilde{H},
\end{equation}
where $H_{ss}$ is a function such that $X_{H_{ss}}$ is the semisimple part
of the linear part of $X_H$, and
\begin{equation}
\{H,H_{ss}\} = 0 .
\end{equation}
\end{thm}

{\it Proof}. For any function $f$ on $\K^{2l+m}$, we can write
$X_f=X_f^{symp}+X_f^{\fg}$ where $X_f^{symp}$ (resp. $X_f^{\fg}$) is the
Hamiltonian vector field of $f$ with respect to $\Pi_{symp}$ (resp.
$\Pi_{\fg}$). We can write $H=\sum_{p,q}H^{p,q}$ where $H^{p,q}$ is a
polynomial of degree $p$ in $x,y$ and of degree $q$ in $z$.

A difficulty of our situation comes from the fact that $\Pi$ is not
homogeneous. If $p>0$ then $X_{H^{p,q}}$ is not a homogeneous vector field
but the sum of a homogeneous vector field of degree $p+q$ (given by
$X^\fg_{H^{p,q}}$) and a homogeneous vector field of degree $p+q-1$ (given
by $X^{symp}_{H^{p,q}}$). Note that $X^\fg_{H^{0,q}}$ is homogeneous of
degree $q$ and of course $X^{symp}_{H^{0,q}}=0$.

Denoting by $X^{(1)}$ the linear part of $X_H$, we have
\begin{equation}
X^{(1)}=X_{H^{0,1}}+X_{H^{2,0}}+X_{H^{1,1}}^{symp}\,.
\end{equation}
This linear vector field $X^{(1)}$ is \emph{not} a Hamiltonian vector
field in general, but we will show that its semisimple part is
Hamiltonian.

By complexifying the system if necessary, we will assume that $\K = \C$.
By a linear canonical change of coordinates, we can suppose that the
semisimple part of $X_{H^{2,0}}$ is $X_{h_2}$ where $h_2(x,y)=\sum_{j=1}^l
\gamma_i x_jy_j$ ($\gamma_j \in \C$) and that the semisimple part of
$X_{H^{0,1}}$ is $X_{h_1}$ where $h_1$ belongs to a Cartan subalgebra
$\fh$ of $\fg$. We write :
\begin{equation}
X_{h_2}=-\sum_{j=1}^l \gamma_jx_j\frac{\partial }{\partial x_j}+
\sum_{j=1}^l \gamma_jy_j\frac{\partial }{\partial y_j}\quad
\mathrm{and}\quad X_{h_1}=\sum_{j=1}^m \alpha_j z_j\frac{\partial
}{\partial z_j}\,. \label{eqn:sspart}
\end{equation}
Remark that we can assume that $\alpha_{s+1}=\hdots=\alpha_m=0$ where
$m-s$ is the dimension of the Cartan subalgebra $\fh$. Denote
$\alpha=(\alpha_1,\hdots,\alpha_m)$ and
$\gamma=(\gamma_1,\hdots,\gamma_l)$. If $\lambda, \mu \in \Z_+^l$ and $\nu
\in \Z_+^m$ then
\begin{equation}
\{h_1+h_2,x^\lambda y^\mu z^\nu\}= ( \langle \gamma , \mu-\lambda \rangle
+ \langle \alpha, \nu \rangle)\, x^\lambda y^\mu z^\nu\,,
\label{eqn:hamilt1}
\end{equation}
where, for example,  $\langle \alpha, \nu \rangle = \sum\alpha_j\nu_j$
denotes the standard scalar product of $\alpha$ and $\nu$. In particular,
$\{h_1 + h_2, . \}$ acts in a ``diagonal'' way on monomials.

We can arrange so that, written as a matrix, the terms coming from
$X_{H^{0,1}- h_1}$, $X_{H^{2,0}-h_2}$ and $X^{symp}_{H^{1,1}}$ in the
expression of $X^{(1)}$ are off-diagonal upper-triangular (and the terms
coming from $X_{h_1 + h_2}$ are on the diagonal).

If $\{h_1+h_2,H^{1,1}\} = 0$, then $[X_{h_1+h_2},X_{H^{1,1}}] = 0$, and
$[X_{h_1+h_2},X^{symp}_{H^{1,1}}] = 0$ because $X_{h_1+h_2}$ is linear and
$X^{symp}_{H^{1,1}} = 0$ is the linear part of $X_{H^{1,1}}$, and, as a
consequence, $X_{h_1+h_2}$ is the semisimple part of $X^{(1)}$.

If $\{h_1+h_2,H^{1,1}\} \neq 0$ then we can apply some canonical changes
of coordinates to make (the new) $H^{1,1}$ commute with $h_1+h_2$ as
follows. According to (\ref{eqn:hamilt1}), there exist two polynomials
$G^{1,1}_{(1)}$ and ${\widetilde{G^{1,1}_{(1)}}}$ of degree 1 in $x,y$ and
1 in $z$ such that
\begin{equation}
\{h_1+h_2,G^{1,1}_{(1)}\}+{\widetilde{G^{1,1}_{(1)}}}=H^{1,1}
\label{eqn:hamilt2}
\end{equation}
and
\begin{equation} \label{eqn:hamilt3}
\{h_1+h_2,{\widetilde{G^{1,1}_{(1)}}} \} = 0 \ .
\end{equation}

Remark that, for any homogeneous polynomials $K^{0,1}, K^{2,0}, K^{1,1}$
of corresponding degrees in $(x,y)$ and $z$, we have
\begin{equation} \label{eqn:doublebracket0}
[X_{K^{1,1}}^{symp}, [X_{K^{1,1}}^{symp},X_{K^{0,1}}]] =
[X_{K^{1,1}}^{symp}, [X_{K^{1,1}}^{symp},X_{K^{2,0}}]] = 0 \ .
\end{equation}

Denote by $F_1 = H^{0,1} + H^{2,0} - h_1 - h_2$ the ``nilpotent part'' of
$H^{0,1} + H^{2,0}$.

Change the coordinate system by the push-forward of the time-1 flow
$\varphi_{(1)} = \exp X_{G^{1,1}_{(1)}}$ of the Hamiltonian vector field
$X_{G^{1,1}_{(1)}}$, i.e., $x^{new}_i = x_i \circ \varphi_{(1)}$ and so
on. The new coordinate system is still a canonical coordinate system,
because $\varphi_{(1)}$ preserves the Poisson structure $\Pi$. By this
canonical change of coordinates, we can replace $H$ by $H \circ
\varphi_{(1)}$, and $X_H$ by
\begin{equation} \label{eqn:XHnew}
X_H^{new} = {\varphi_{(1)}}_\ast X_H= X_H + [X_{G^{1,1}_{(1)}},X_H] +
{1\over 2} [X_{G^{1,1}_{(1)}}, [X_{G^{1,1}_{(1)}},X_H]] + \hdots
\end{equation}
It follows from (\ref{eqn:XHnew}) and (\ref{eqn:doublebracket0}) that the
linear part of $X_H^{new}$ is
\begin{equation}
X_{h_1+h_2}+X_{F_1}+X^{symp}_{{\widetilde{G^{1,1}_{(1)}}}}+
X^{symp}_{\{G^{1,1}_{(1)},F_1\}} .
\end{equation}
In particular, by the above canonical change of coordinates, we have
replaced $H^{1,1} = \{h_1+h_2,G^{1,1}_{(1)}\}+{\widetilde{G^{1,1}_{(1)}}}$
by $\{G^{1,1}_{(1)},F_1\}+{\widetilde{G^{1,1}_{(1)}}}$, while keeping
$h_1,h_2$ and $F_1$ intact. (Note that $\{G^{1,1}_{(1)},F_1\}$ is
homogeneous of degree 1 in $(x,y)$ and degree 1 in $z$).

By (\ref{eqn:hamilt1}), (\ref{eqn:hamilt2}), (\ref{eqn:hamilt3}) we can
write
\begin{equation}
G^{1,1}_{(1)}=\{h_1+h_2,G^{1,1}_{(2)}\},
\end{equation}
which, together with $\{h_1+h_2,F_1\}=0$, gives
\begin{equation}
\{G^{1,1}_{(1)},F_1\}=\{h_1+h_2,\{G^{1,1}_{(2)},F_1\}\} .
\end{equation}
In other words, the new $H^{1,1}$ is $\{h_1+h_2,\{G^{1,1}_{(2)},F_1\}\}
+{\widetilde{G^{1,1}_{(1)}}}$.

Repeating the above process, with the help of the time-1 flow
$\varphi_{(2)}=\exp X_{\{G^{1,1}_{(2)},F_1\}}$ of the Hamiltonian vector
field of ${\{G^{1,1}_{(2)},F_1\}}$, we can replace $H^{1,1}$ by
\begin{equation}
\{h_1+h_2,\{{\{G^{1,1}_{(3)},F_1\},F_1}\}\} +{\widetilde{G^{1,1}_{(1)}}},
\end{equation}
and so on.

Since ${F_1}$ is ``nilpotent'', by iterating the above process a finite
number of times, we can replace $H^{1,1}$ by
${\widetilde{G^{1,1}_{(1)}}}$, i.e. make it commute with $h_1 + h_2$. So
we can assume that $\{h_1+h_2,H^{1,1}\} = 0$. Then
\begin{equation}\label{eqn:sscomp}
H_{ss}=h_1+h_2
\end{equation}
is a function such that $X_{H_{ss}}$ is the semisimple part of the linear
part of $X_{H}$.

Now let us deal with higher degree terms. Write
\begin{equation}
X_H=X_{H_{1}} + X_{H_{2}} + X_{H_{3}} + \hdots ,
\end{equation}
where each $H_k$ is of the type
\begin{equation}
H_{k} = H^{0,k} + \sum_{p \geq 1} H^{p,k+1-p} .
\end{equation}
(For example, $H_1 = H^{0,1} + H^{2,0} + H^{1,1} = H_{ss} + F_1 +
H^{1,1}$).

By recurrence, assume that, for some $r \geq 2$, we have $\{H_{ss},H_{k}\}
= 0$ for all $k \leq r-1$. We will change $H_r$ by a canonical coordinate
transformation to get the same equality for $k = r$.

 In order to put $H_r$ in normal form, we use the same method that we used
 to normalize $H^{1,1}$. Similarly to (\ref{eqn:hamilt2}), we can write
\begin{equation}
H_r = \{H_{ss},K_r\} + \tilde{K}_r ,
\end{equation}
where $K_r$ and $\tilde{K}_r$ are of the same type as $H_r$ (i.e., they
are sums of monomials of bidegrees $(0,r)$ and $(p,r+1-p)$ with $p > 0$),
$\{H_{ss},\tilde{K}_r\} = 0$. Note $K_r$ can be written as $K_r =
\{H_{ss},K_{(2)r}\}$ for some $K_{(2)r}$.

The canonical coordinate transformation given by the time-1 flow $\exp
X_{K_r}$ of $X_{K_r}$ leaves $H_1,\hdots, H_{r-1}$ intact, and changes
$H_r =  \tilde{K}_r + \{H_{ss},K_r\}$ to the sum of $\tilde{K}_r$ with
the terms of appropriate bidegrees in $\{K_r,F_1+H^{1,1}\}$. We will write
it as
\begin{equation}
\tilde{K}_r + \{K_r,F_1+H^{1,1}\} \ mod \ (terms \ of \ higher \
bidegrees) .
\end{equation}
It can also be written as
\begin{equation}
\tilde{K}_r + \{ H_{ss},\{K_{(2)r},F_1+H^{1,1}\}\} \ mod \ (terms \ of \
higher \ bidegrees) .
\end{equation}
Now apply the canonical coordinate transformation given by $\exp
X_{\{K_{(2)r},F_1+H^{1,1}\}}$, and so on. Since $F_1+H^{1,1}$ is
``nilpotent'', after a finite number of coordinate transformations
like that, we can change $H_r$ to $\tilde{K}_r$, which commutes
with $H_{ss}$. Denote the composition of these coordinate changes
(for a given $r$) as $\phi_r$. Note that $\phi_r$ is of the type
\begin{equation}
\phi_r=Id+\, terms \,of\, degree \, \geq r
\end{equation}
Thus, the sequence of local or formal Poisson-structure-preserving
diffeomorphisms $(\Phi_r)_{r \geq 2}$, where $\Phi_r = \phi_r \circ \hdots
\circ \phi_2$, converges formally and gives a formal normalization of $H$.

Finally, notice that, in the real case ($\K = \R)$, by an argument similar
to the one given in the previous section, all canonical coordinate
transformations constructed above can be chosen real.

Theorem \ref{thm:formal-Birkhoff} is proved. \QED

\begin{rem} In Theorem \ref{thm:formal-Birkhoff}, if we forget the Lie algebra
$\fg$ and just keep the symplectic structure, then we recover the
classical Birkhoff normalization for Hamiltonian vector fields on
symplectic manifolds (see, e.g.,
\cite{Birkhoff,Bruno2,Siegel-Moser,Zung-Birkhoff}). On the other hand, if
we forget the symplectic part and just deal with $\fg^\ast$ then we get
the following result as a particular case:

\begin{cor} Let $h$ be a local analytic or formal function , with $h(0)=0$ and $dh(0)\neq
0$, on the dual $\fg^\ast$ of a semisimple Lie algebra with the associated
Lie-Poisson structure. Then the Hamiltonian vector field $X_h$ admits a
formal Poincaré-Birkhoff normalization, i.e., there exists a formal
coordinate system in which the Poisson structure is linear and in which we
have
$$
\{h,h_{ss}\} =0\,,
$$
where $h_{ss}$ is the semisimple part of $dh(0)$ in $\fg$.
\end{cor}

\end{rem}

\begin{example}
The monomials $x^\lambda y^\mu z^\nu$ such that $\langle
\gamma,\mu-\lambda \rangle+ \langle \alpha , \nu \rangle =0$ in
(\ref{eqn:hamilt1}) may be called {\it resonant} terms. In the two following examples
we give the set of all resonant terms in the case of a trivial symplectic part.

{\it a)} $\fg=sl(2)$. In this case, a Cartan subalgebra $\mathfrak{h}$ of
$\fg$ is of dimension 1 and there are only two roots $\{\alpha,-\alpha\}$.
Denote by ${z_1,z_2,z_3}$ a basis of $\fg$ (or a coordinate system on
$\fg^\ast$) such that $z_1$ (resp. $z_2$) spans the root space associated
to $\alpha$ (resp. $-\alpha$) and $z_3$ spans the Cartan subalgebra. We
suppose that in the decomposition (\ref{eqn:sscomp}) we have $h_1=z_3$.
Then the resonant terms are formal power expansion
in the variables $\omega=z_1z_2$ and $z_3$.

{\it b)} $\fg=sl(3)$. Here a Cartan subalgebra $\mathfrak{h}$ is of
dimension 2 (see for instance \cite{Knapp}). There are 6 roots
$\{\alpha_1,\alpha_2,\alpha_3,-\alpha_1,-\alpha_2,-\alpha_3\}$ and the
relations between these roots are of type
\begin{equation}
\sum_i a_i\alpha_i - \sum_i b_i\alpha_i=0 \,,
\end{equation}
with
\begin{equation}
a_1-b_1=a_2-b_2=a_3-b_3\,.
\end{equation}
If $\{\xi_1,\xi_2,\xi_3,\zeta_1,\zeta_2,\zeta_3,z_1,z_2\}$ is a basis of
$\fg$ such that $\xi_j$ (resp. $\zeta_j$) spans the root space associated
to $\alpha_j$ (resp. $-\alpha_j$) and $\{z_1,z_2\}$ spans $\mathfrak{h}$,
then supposing that in the decomposition (\ref{eqn:sscomp}) $h_1$ is a
linear combination of $z_1$ and $z_2$ we may write the  resonant terms as
formal power expansion formed by monomials of type
\begin{equation}
\xi_1^{a_1}\zeta_1^{b_1} \xi_2^{a_2}\zeta_2^{b_2}
\xi_3^{a_3}\zeta_3^{b_3} z_1^{u_1}z_2^{u_2}
\end{equation}
with $a_1-b_1=a_2-b_2=a_3-b_3$.
\end{example}

\subsection{Analytic normalization for integrable Hamiltonian systems}
Here, we assume that we work in the {\it complex analytic setup}.

Recall that we wrote in (\ref{eqn:sspart}),
\begin{equation}
X_{h_2}=-\sum_{j=1}^l \gamma_jx_j\frac{\partial }{\partial x_j}+
\sum_{j=1}^l \gamma_jy_j\frac{\partial }{\partial y_j}\quad
\mathrm{and}\quad X_{h_1}=\sum_{j=1}^m \alpha_j z_j\frac{\partial
}{\partial z_j}\,,
\end{equation}
and we had put $\alpha=(\alpha_1,\hdots,\alpha_m)\in\K^m$ and
$\gamma=(\gamma_1,\hdots,\gamma_l)\in\K^l$.

Let $\mathcal{R}\subset \Z^{2l+m}$ be the sublattice of
$\Z^{2l+m}$ formed by the vector $u\in\Z^{2l+m}$ written as
$u=(\lambda,\mu,\nu)$, with $\lambda$ and $\mu$ in $\Z^l$ and
$\nu$ in $\Z^m$, and such that
\begin{equation}
\langle(-\gamma,\gamma,\alpha)\, , \,(\lambda,\mu,\nu)\rangle= -\sum \gamma_j
\lambda_j + \sum \gamma_j \mu_j + \sum \alpha_j \nu_j =0\,.
\end{equation}
Of course, the elements $(\lambda,\mu,\nu)$ of $\mathcal{R}$
correspond to the resonant monomials i.e. terms of type $x^\lambda
y^\mu z^\nu$ such that $\{H_{ss},x^\lambda y^\mu z^\nu\}=0$. The
dimension of $\mathcal{R}$ may be called the {\it degree of
resonance} of $H$.

Now, we consider the sublattice $\mathcal{Q}$ of $\Z^{2l+m}$
formed by vectors $a\in \Z^{2l+m}$ such that $\langle a\,|\,u \rangle =0$ for all
$u$ in $\mathcal{R}$. Let $\{\rho^{(1)},\hdots,\rho^{(r)}\}$ be a
basis of $\mathcal{Q}$. The dimension $r$ of $\mathcal{Q}$ is
called the {\it toric degree} of $X_H$ at 0. We then put for all
$k=1,\hdots,r$
\begin{equation}
Z_k=\sum_{j=1}^l \rho^{(k)}_j x_j \frac{\partial }{\partial x_j}+
\sum_{j=1}^l \rho^{(k)}_{l+j} y_j \frac{\partial }{\partial y_j}+
\sum_{j=1}^m \rho^{(k)}_{2l+j} z_j \frac{\partial }{\partial
z_j}\,. \label{eqn:deftore}
\end{equation}

The vector fields $iZ_1,\hdots,iZ_r$ ($i=\sqrt{-1}$) are periodic
with a real period in the sense that the real part of these vector
fields is a periodic real vector field in
$\C^{2l+m}=\R^{2(2l+m)}$ ; they commute
pairwise and are linearly independent almost everywhere. Moreover,
the vector field $X_{H_{ss}}$ is a linear combination (with
coefficients in $\C$ a priori) of the $iZ_k$. We also have the
trivial following property

\begin{lem}
If $\Lambda$ is a $p$-vector ($p\geq 0$) then we have the
equivalence
$$
[X_{H_{ss}},\Lambda]=0 \Leftrightarrow [Z_k,\Lambda]=0\; \forall
k=1,\hdots,r
$$
\label{lem:commutZk}
\end{lem}

{\it Proof :} We just give here the idea of the proof of this
lemma supposing that $\Lambda$ is a 2-vector for instance ; but it
works exactly in the same way for other multivectors. If $Y$ is a
vector field of type $\sum_{j=1}^l a_j x_j \frac{\partial
}{\partial x_j}+ \sum_{j=1}^l a_{l+j} y_j \frac{\partial
}{\partial y_j}+ \sum_{j=1}^m a_{2l+j} z_j \frac{\partial
}{\partial z_j}$ and $\Lambda$ of type $\Lambda=x^\lambda y^\mu
z^\nu \frac{\partial }{\partial x_u}\wedge \frac{\partial
}{\partial x_v}$, then
\begin{equation}
[Y,\Lambda]= \langle a\, , \, (\lambda,\mu,\nu)-(1_u,1_v,0) \rangle \Lambda
\end{equation}
where $1_u=(0,\hdots,1,\hdots,0)$ is the vector of $\Z^l$ whose
unique nonzero component is the $u$-component. Of course we get
the same type of relation with 2-vectors in $\frac{\partial
}{\partial x}\wedge \frac{\partial }{\partial z}$, $\frac{\partial
}{\partial x}\wedge \frac{\partial }{\partial x}$, etc.... Using
this remark and the definition of the vectors
$\rho^{(1)},\hdots,\rho^{(r)}$, the equivalence of the lemma is
direct. \QED

According to this Lemma, since $X_{H_{ss}}$ preserves the Poisson
structure, the vector fields $Z_1,\hdots,Z_r$ will be Poisson
vector fields for
$(\C^{2l}\times\C^m,\{\,,\,\}_{symp}+\{\,,\,\}_{\fg^\ast})$. But
according to Proposition \ref{prop:kunneth} (see the Appendix),
the Poisson cohomology space
$H^1(\C^{2l}\times\C^m,\{\,,\,\}_{symp} +\{\,,\,\}_{\fg})$ is
trivial therefore, these vector fields are actually Hamiltonian :
\begin{equation}
Z_k=X_{F_k}\quad \forall k=1,\hdots,r\,.
\end{equation}

Finally, we have $r$ periodic Hamiltonian linear vector fields
$iZ_k$ which commute pairwise, are linearly independent almost
everywhere. The real parts of these
 vector fields generate a Hamiltonian action of the real torus $\mathbb{T}^r$
on $(\C^{2n}\times\C^m,\{\,,\,\}_{symp} + \{\,,\,\}_\fg)$.
With all these notations, we can state the following proposition :

\begin{prop}
With the above notation, the following conditions are equivalent~:

a) There exists a holomorphic Poincar\'e-Birkhoff normalization of
$X_H$ in a neighborhood of 0 in $\C^{2l+m}$.

b) There exists an analytic Hamiltonian action of the real torus
$\mathbb{T}^r$ in a neighborhood of 0 in $\C^{2l+m}$, which
preserves $X_H$ and whose linear part is generated by the
(Hamiltonian) vector fields $iZ_k$, $k=1,\hdots,r$.
\label{prop:equivPB-Torus}
\end{prop}

{\it Proof :} Suppose that $H$ is in holomorphic
Poincar\'e-Birkhoff normal form. By Lemma \ref{lem:commutZk},
since $\{H,H_{ss}\}=0$, the vector fields $iZ_k$ preserve $X_H$.

Conversely, if the point {\it b)} is satisfied, then according to the
 holomorphic version of the Splitting Theorem (see \cite{Eva-Zung}) we can
consider that the action of the torus is ``diagonal'', i.e. the
product of an action on $(\C^{2l},\{\,,\}_{symp})$ by an action on
$(\C^m,\{\,,\,\}_\fg)$ and moreover that the action on the
symplectic part is linear. According to Proposition
\ref{prop:linearaction} (in Appendix), we can linearize the second
part of the action by a Poisson diffeomorphism. We then can
consider that the action of $\mathbb{T}^r$ is generated by the
vector fields $iZ_k$, $k=1,\hdots,r$. This action preserves $X_H$
then we have $[iZ_k,X_H]=0$ for all $k$. To conclude, just recall
that $X_{H_{ss}}$ is a linear combination of the $Z_k$. \QED

Now, we are going to use Proposition \ref{prop:equivPB-Torus} to
clarify a link between the integrability of a Hamiltonian vector
field $X_H$ on an analytic Poisson manifold $(\K^n,\{\,,\,\})$ and
the existence of a convergent Poincar\'e-Birkhoff normalization.
Recall first the definition (see for instance \cite{Zung-survey})
of the word {\it integrability} used here :

\begin{defn}
A Hamiltonian vector field $X_H$ on a Poisson manifold $(M,\Pi)$
(of dimension $n$) is called {\it integrable} (in the generalized
Liouville sense) if there exist $p$ ($1\leq p\leq n$) Hamiltonian
vector fields $X_1=X_H,X_2,\hdots,X_p$ and $n-p$ functions
$f_1,\hdots,f_{n-p}$ such that

{\it a)} The vector fields commute pairwise, i.e.
\begin{equation}
[X_i,X_j]=0\;\forall i,j=1,\hdots,p\,,
\end{equation}
and they are linearly independent almost everywhere, i.e.
\begin{equation}
X_1\wedge\hdots\wedge X_p\neq 0\,.
\end{equation}
{\it b)} The functions are common first integrals for
$X_1,\hdots,X_p$ :
\begin{equation}
X_i(f_j)=0\;\forall i,j\,,
\end{equation}
and they are functionally independent almost everywhere :
\begin{equation}
df_1\wedge\hdots\wedge df_{n-p}\neq 0\,.
\end{equation}
\end{defn}

Of course this definition has a sense in the smooth category as
well as in the analytic category. We can speak about {\it smooth}
or {\it analytic} integrability.

\begin{thm}
Any analytically integrable Hamiltonian vector field in a
neighborhood of a singularity on an analytic Poisson manifold
admits a convergent Poincar\'e-Birkhoff normalization
\label{thm:integrability-cvgtPB}
\end{thm}

{\it Proof :} We can assume (see the beginning of the section)
that we work in the neighborhood of 0 in
$$
(\C^{2l+m},\{\,,\,\})=(\C^{2l},\{\,,\,\}_{symp})\times
(\fg^\ast,\{\,,\,\}_\fg)
$$
where $\{\,,\,\}_{symp}$ is a symplectic Poisson structure and
$\fg$ is a semisimple Lie algebra and $\{\,,\,\}_\fg$ the standard
Lie-Poisson structure on $\fg ^\ast$. If $X_H$ is integrable then,
forgetting one moment the Hamiltonian feature, Theorem 1.1 and
Proposition 2.1 in \cite{Zung-Dulac} give the existence of an
action of a real torus $\mathbb{T}^r$ on $(\K^{2l+m},0)$ generated
by vector fields $Y_1,\hdots,Y_r$ ($r$ is the toric degree of
$X_H$) where the linear parts of these vector fields are the
$iZ_k$ (see \ref{eqn:deftore}), and which preserves $X_H$.
Moreover, the semisimple part $X^{ss}_H$ of $X_H$ is a linear
combination of the $Y_j$ : $X^{ss}_H=\sum_j \beta_j Y_j$ without any
resonance relation between the $\beta_j$. Now, let us recall that
we work in a Poisson manifold with a Hamiltonian vector field.
Since the vector field $X_H$ preserves the Poisson structure, its
semisimple part also does and then we will have $[Y_j,\Pi]=0$ for
all $j=1,\hdots,r$. Therefore, the action of the torus also
preserves the Poisson structure. Proposition
\ref{prop:equivPB-Torus} allows to conclude.\QED

\begin{rem}
If we suppose that $H$ and the Poisson structure are real then it
is natural to ask if all that we made is still valid. Note that in
this case, we can consider $H$ (and the Poisson structure) as
complex analytic, with real coefficients.

Actually, in the same way as in \cite{Zung-Dulac,Zung-Birkhoff},
we conjecture that we have the equivalence:

{\it A real analytic Hamiltonian vector field $X_H$ with respect
to a real analytic Poisson structure admits a local real analytic
Poincar\'e-Birkhoff normalization iff it admits a local
holomorphic Poincar\'e-Birkhoff normalization.}
\end{rem}

\section{Appendix}
In this appendix, we give a proof of auxiliary results used in the
previous sections. We first compute the first Poisson cohomology space of
the Poisson manifold we consider in Section \ref{section:Hamiltonian}.
Suppose that $\Pi_S$ is a symplectic (i.e. nondegenerate) Poisson
structure on $\K^{2l}$ ($\K$ is $\R$ or $\C$). If
$(x_1,\hdots,x_l,y_1,\hdots,y_l)$ are coordinates on $\K^{2l}$, we can
write
$$\Pi_S=\sum_{i=1}^l \frac{\partial }{\partial x_i}\wedge
\frac{\partial }{\partial y_i}\,.
$$
Let $\fg$ be a $m$-dimensional (real or complex) semisimple Lie
algebra and consider $\Pi_{\fg}$ the corresponding linear Poisson
structure on $\K^m$. Suppose that $(z_1,\hdots,z_m)$ are
coordinates on $\K^m$. We then show the following :

\begin{prop}
Under the hypotheses above, if
$\mathrm{H}^1(\K^{2l}\times\K^m,\Pi_S+\Pi_{\fg})$ denotes the first (formal
or analytic) Poisson cohomology space of the product of
$(\K^{2l},\Pi_S)$ by $(\K^m,\Pi_{\fg})$ then
$$
\mathrm{H}^1(\K^{2l}\times\K^m,\Pi_S+\Pi_{\fg})=\{0\}\,.
$$
\label{prop:kunneth}
\end{prop}

{\it Proof :} If $X$ is a (formal or analytic) vector field on
$\K^{2l}\times\K^m$ we write $X=X^S+X^{\fg}$ where $X^S$ is a
vector field which only has components in the $\frac{\partial
}{\partial x_i}$ and $\frac{\partial }{\partial y_i}$ and, in the
same way, $X^{\fg}$ only has components in the $\frac{\partial
}{\partial z_i}$. Before computing the Poisson cohomology space,
let us make the following two remarks :

If $[X^S,\Pi_S] = 0$ then $X^S=[f,\Pi_S]$ where $f$ is a (formal
or analytic) function on $\K^{2l+m}$. Indeed, recalling that
(because $\Pi_S$ is symplectic) the Poisson cohomology of
$(\K^{2l},\Pi_S)$ is isomorphic to the de Rham cohomology of
$\K^{2l}$ (see for instance \cite{Vaisman}), the relation
$[X^S,\Pi_S] = 0$ may be translated as $d\alpha=0$ where $\alpha$
is a 1-form on $\K^{2l}$ depending (formally or analytically) on
parameters $z_1,\hdots,z_m$. Then we can write $\alpha=d f$ where
$f$ is a function on $\K^{2l}$ depending (formally or
analytically) on parameters $z_1,\hdots,z_m$.

In the same way, if $[X^{\fg},\Pi_S] = 0$ then, writing
$X^{\fg}=\sum_i X^{\fg}_i(x,y,z) \frac{\partial }{\partial z_i}$,
we get $[X^{\fg}_i,\Pi_S]=0$ for all $i$. Thus, each $X^{\fg}_i$
depends only on $z$. Indeed, here $X^{\fg}_i$ may be seen as a
function on $\K^{2l}$ depending (formally or analytically) on
parameters $z_1,\hdots,z_m$ such that $d X^{\fg}_i=0$.\\

Now if $X=X^S+X^{\fg}$ is a vector field on $\K^{2l}\times\K^m$,
it is easy to see that $[X,\Pi_S+\Pi_{\fg}]=0$ is equivalent to the
three equations

\begin{eqnarray}
0 &=& [X^S,\Pi_S]  \label{eqn:SS}\\
0 &=& \left[X^S,\Pi_{\fg}\right]+\left[X^{\fg},\Pi_S\right] \label{eqn:SG} \\
0 &=& \left[ X^{\fg},\Pi_{\fg} \right] \label{eqn:GG}
\end{eqnarray}

According to the first remark we made above, equation
(\ref{eqn:SS}) gives $X^S=[f,\Pi_S]$ where $f$ is a (formal or
analytic) function on $\K^{2l+m}$. Now, replacing $X^S$ by
$[f,\Pi_S]$ in (\ref{eqn:SG}) and using the graded Jacobi identity
of the Schouten bracket, we get
\begin{equation}
\big[ X^{\fg}-[f,\Pi_{\fg}], \Pi_S\big] =0\,.
\end{equation}
Since $X^{\fg}-[f,\Pi_{\fg}]$ is a vector field which only has
component in $\frac{\partial}{\partial z}$, the second remark we
made above gives
\begin{equation}
X^{\fg} = [f,\Pi_{\fg}] +Y
\end{equation}
where $Y$ is a vector field on $\K^m$ (i.e. only has components in
$\frac{\partial}{\partial z}$ and whose coefficients are functions
of $z$). Finally, (\ref{eqn:GG}) gives $[Y,\Pi_{\fg}]=0$ i.e. $Y$
is a 1-cocycle for the Poisson cohomology of $(\K^m,\Pi_{\fg})$.
Since the Lie algebra $\fg$ is semisimple, the Poisson cohomology
space $\mathrm{H}^1(\K^m,\Pi_{\fg})$ is trivial (see for instance
\cite{Conn-Analytic1984}). We then obtain $Y=[h,\Pi_{\fg}]$ where
$h$ is a function on $\K^m$.

To resume, we get
\begin{equation}
X=X^S+X^{\fg}=[f+h,\Pi_S+\Pi_{\fg}]\,,
\end{equation}
which means that $X$ is a 1-cobord for the Poisson cohomology of
$(\K^{2l}\times\K^m,\Pi_S+\Pi_{\fg})$.\QED\\

The second result is an analytic version of a smooth linearization
theorem due to V. Ginzburg. In the Appendix of  \cite{Ginzburg},
he states that the $G$-action of a compact Lie group on a Poisson
manifold $(P,\Pi)$ (everything is smooth here), fixing a point $x$
of $P$ and such that the Poisson structure is linearizeable at
$x$, can be linearized by a diffeomorphism which preserves the
Poisson structure. Here, we state the following~:

\begin{prop}
Consider an analytic action of a compact (analytic) Lie group on
$(\K^n,\Pi)$ ($\K$ is $\R$ or $\C$), where $\Pi$ is an analytic
Poisson structure on $\K^n$. Suppose that the action fix the
origin 0 and that the Poisson structure is linearizable at 0.
Then, the action can be linearized by a Poisson diffeomorphim.
\label{prop:linearaction}
\end{prop}

{\it Proof :} The proof is the same as in the smooth case : we use
Moser's path method. If $g$ is an element of $G$, we put
$\varphi^g$ the corresponding diffeomorphism of $\K^n$ and
$\varphi^g_{lin}$ its linear part at 0. We construct a path of
analytic actions of $G$ on $(\K^n,\Pi)$ given by the following
diffeomorphisms :
$$
\varphi_t^g(x)= \left\{
\begin{array}{ccc}
\varphi^g(tx)/t & \mathrm{if} & 0<t\leq 1\\
\varphi^g_{lin}(x) & \mathrm{if} & t=0
\end{array} \right.
$$
for any $g$ in $G$ and $x$ in $\K^n$. These actions preserve $\Pi$
and fix 0. We want now to show that there exists a path of
diffeomorphisms $\psi_t$, with $\psi_0=Id$,  preserving the
Poisson structure $\Pi$ and such that
\begin{equation}
\psi_t\circ\varphi_t^g\circ\psi_t^{-1}=\varphi_0^g=\varphi_{lin}^g\,,
\label{eqn:Moser-action}
\end{equation}
for all $t\in [0,1]$ and all $g$ in $G$.

Let $C_t(g)$ be the time-depending vector field associated to
$\varphi_t^g$ :
\begin{equation}
C_t(g)(\varphi_t^{g}(x))=\frac{\partial \varphi_t^g}{\partial t}
(x)\,.
\end{equation}

Derivating the condition (\ref{eqn:Moser-action}), we are led to
look for a time-depending vector field $X_t$ (corresponding to
$\psi_t$) verifying
\begin{equation}
C_t(g)={\varphi_t^g}_\ast X_t - X_t\,,
\label{eqn:cohomologicalcond}
\end{equation}
for all $t\in[0,1]$ and all $g$ in $G$.

We put
\begin{equation}
X_t=-\int_G {\varphi_t^h}_\ast C_t(h) dh \,,
\end{equation}
$dh$ is a bi-invariant Haar measure on $G$ such that the volume of
$G$ is 1. This vector field is analytic and depends smoothly on
$t$. Moreover, since each $C_t(h)$ preserves the Poisson structure
$\Pi$, so does $X_t$. Finally, one can check that $X_t$ satisfies
the condition (\ref{eqn:cohomologicalcond}).\QED

{\it Proof of Proposition \ref{prop:fullmeasure}}.  We denote by $\alpha$ the linear
application from the Cartan subalgebra $\fh$ to $\K^n$ defined by
$\alpha(h)=(\alpha_1(h),\hdots,\alpha_n(h))$ for any $h$ in $\fh$ and by
$W$ its image. We show that the subset of $W$
formed by the elements $\gamma$ such that the $\omega_d(\gamma)$
(defined as in (\ref{eqn:omega_d}) replacing
$\langle \alpha_i \,,\,h_1 \rangle  $ by $\gamma_i$) do not
satisfy the $\omega$-condition is of measure 0 (in $W$). Since $\alpha$ is a linear
surjection from $\fh$ to $W$, it will show Proposition \ref{prop:fullmeasure}.

Note that if $\gamma\in\K^n$ satisfies the condition (which is a
condition of type "Siegel")
\begin{equation}
(\exists c>0)\;(\forall \lambda\in \Z_+^n),\; s.t.\;
\big|  |\lambda| - 1 + \langle \gamma \,,\, \lambda \rangle \big| \geq \frac{c}{|\lambda|^s}\,,
\label{eqn:Siegel}
\end{equation}
where $s>n$, then  $\omega_d(\gamma)$
satisfies the $\omega$-condition
(\ref{eqn:Bruno}). We then show that the set of the $\gamma$ in $W$ which
do not satisfy the condition (\ref{eqn:Siegel}) is of measure 0 in $W$.

For any positive integer $k$ and any positive real number $c$, if $\|\,\|$
denotes the norm associated to $ \langle \,,\, \rangle$, we put
\begin{eqnarray*}
W_k &=& \big\{\gamma\in W\;\big|\; \|\gamma\|\leq k \big\} \\
V_{c} &=& \big\{ \gamma\in\K^n \, \big| \, (\exists \lambda \in
\Z_+^n)\; s.t.\;
\big|  |\lambda| - 1 + \langle \gamma \,,\, \lambda \rangle \big|
\, \leq \frac{c}{ |\lambda|^s } \, \big\} \\
V &=& \cap_{c>0} V_{c}
\end{eqnarray*}
Actually, we show here that $W_1\cap V$ is of measure 0 but the same
technic works to prove that $W_{k}\cap V$ is also of measure 0 for each $k$.
Therefore $\cup_k (W_k\cap V)$ is of measure 0 too, which proves the proposition.

Now, for any $\lambda$ in $\Z_+^n$ we consider the affine
subspace $\mathcal{V}_\lambda$ of $\K^n$
formed by the vectors $\gamma$ such that
$ \langle \gamma \,,\, \lambda \rangle = 1-|\lambda|$
and we put for $c>0$,
\begin{equation}
\mathcal{V}_{\lambda,c}=\big\{ \gamma\in\K^n\, ;
\, \big|  |\lambda| - 1 + \langle \gamma \,,\, \lambda \rangle \big|
\, \leq \frac{c}{ |\lambda|^s } \, \big\}\,.
\end{equation}
This last set is like a tubular neighborhood of $\mathcal{V}_\lambda$ of
thickness $\frac{2c}{|\lambda|^s}$. We look now at
$K_{\lambda,c}=\mathcal{V}_{\lambda,c}\cap W_1$. If it is not empty, it is a kind of "band" in
$W_1$ of thickness smaller than $S\frac{2c}{|\lambda|^s}$ where $S$ is a
positive constant which only depends on the dimension of $W$ (and on the
metric). Therefore, we get
\begin{equation}
Vol(W_1\cap V_c)\leq \sum_{\lambda\in \Z_+^n}
Vol(K_{\lambda,c}) \leq cS\sum_{\lambda\in \Z_+^n}
\frac{1}{|\lambda|^s}\,.
\end{equation}
This latest sum converges (because $s>n$) and we then get $Vol(W_1\cap
V)=Vol(\cap_{c>0} W_1\cap V_c)=0$.\QED

{\it Proof of Lemma \ref{lem:normes}}. {\it a)} The first inequality of
(\ref{eqn:eqnormes}) is obvious. To prove the second one, we use the
Cauchy inequality
$$
|a_\lambda| \leq \frac{\sup_{z\in D_{\rho}}|f(z)|}{{\rho}^{|\lambda|}}
$$
for all $\lambda$, which induces $|a_\lambda|{\rho^\prime}^{|\lambda|}
\leq \|f\|_{\rho} \big(\frac{\rho^\prime}{\rho}\big)^{|\lambda|}$. The
inequality follows.

The point {\it b)} is obvious.

{\it c)} If $f=\sum_{|\lambda|\geq q} a_\lambda x^\lambda$ then
\begin{eqnarray*}
\Big|\frac{\partial f}{\partial x_j}\Big|_{\rho^\prime} &=&
\sum_{|\lambda|\geq q}
\lambda_j |a_\lambda| {\rho^\prime}^{|\lambda|-1}\\
 &=& \sum_{|\lambda|\geq q}
|a_\lambda| {\rho}^{|\lambda|} \times
\frac{\lambda_j}{\rho^\prime}{\big(\frac{\rho^\prime}{\rho}\big)}^{|\lambda|}
\end{eqnarray*}
When $\rho' = {\Big( \frac{1}{(2d)(2^d)} \Big)}^{1/(2^d+1)} \rho \geq R
> 0$, $q= 2^d+1$ and $d\geq 1$, each number
$\frac{\lambda_j}{\rho^\prime}{\big(\frac{\rho^\prime}{\rho}\big)}^{|\lambda|}$
can be majored by
$$
{2^d+1 \over R} \left( {\Big( \frac{1}{(2d)(2^d)}
\Big)}^{1/(2^d+1)} \right)^{2^d+1}\,.
$$
It is easy to see that these numbers are smaller
than 1, provided that $d$ is large enough. \QED\\

{\it Proof of Lemma \ref{lem:radii}}. {\it a)} Since the sequence
${(r_d)}_d$ decreases and converges to a positive real number $R>0$, we
have $r_d>R$ for all $d$. We write
$r_d-\rho_d=r_d\frac{1}{d^2}>\frac{R}{d^2}$, thus for $d$ sufficiently
large, we get $r_d-\rho_d>\frac{1}{2^d}$.

{\it b)} We have $\rho_d-r_{d+1}=\rho_d\big[1-
{\big(\frac{\omega_{d+1}}{2^{d+1}}\big)}^{\frac{1}{2^{d+1}+1}}\big]$.
Since the sequence ${(\rho_d)}_d$ decreases and converges to $R>0$, we
have $\rho_d>R>0$ for all $d$. We then show that if $d$ is sufficiently
large, then
$$
R \Big[1-
{\big(\frac{\omega_{d+1}}{2^{d+1}}\big)}^{\frac{1}{2^{d+1}+1}}\Big]
>\frac{1}{2^d}\,.
$$
We have
${\big(\frac{\omega_{d+1}}{2^{d+1}}\big)}^{\frac{1}{2^{d+1}+1}}=e^{\gamma_d}$
where $\gamma_d=\frac{1}{2^{d+1}+1} \ln(\frac{\omega_{d+1}}{2^{d+1}})$. By
the $\omega$-condition, the sequence ${(\gamma_d)}_d$ converges to 0 and
is negative for all $d$ sufficiently large. Then, if $\varepsilon$ is a small
positive real number (for instance $\varepsilon=1/2$), we have for all $d$ sufficiently large,
$1-e^{\gamma_d}>-(1-\varepsilon)\gamma_d$. We deduce that
\begin{equation}
R(1-e^{\gamma_d})>-R
(1-\varepsilon)\frac{\ln\big(\frac{\omega_{d+1}}{2^{d+1}}\big)}{2^{d+2}}
\end{equation}
which gives
\begin{equation}
R(1-e^{\gamma_d})>\frac{1}{2^d}\Big[\frac{R(1-\varepsilon)}{4}\big(\ln
(2^{d+1}) - \ln\omega_{d+1}\big)\Big]\,.
\end{equation}
Therefore, for $d$ sufficiently large,
$R(1-e^{\gamma_d})>\frac{1}{2^d}$.\QED \\

\bibliographystyle{amsalpha}

\providecommand{\bysame}{\leavevmode\hboxto3em{\hrulefill}\thinspace}


\end{document}